\documentclass[11pt]{article}
\usepackage{geometry}               
\usepackage[parfill]{parskip}    
\usepackage{amssymb}
\usepackage{amsmath}
\usepackage{color, graphicx}
\usepackage[all]{xy}

\usepackage{calrsfs}
\usepackage[T1]{fontenc} 
\usepackage{textcomp}
\usepackage{times}
 \usepackage[scaled=0.92]{helvet} 
 \usepackage{bbm} 
 
\newcommand{\F}{\mathbf{F}}
\newcommand{\p}{\mathfrak{p}}
\newcommand{\Q}{\mathbf{Q}}
\newcommand{\ab}{{\footnotesize \mathrm{ab}}}
\newcommand{\isomto}{\overset{\sim}{\rightarrow}}
\newcommand{\tr}{\mathrm{tr}}

\begin{document}
{\Large {\textbf{Class field theory as a dynamical system}}}\\[1mm]
{\large \emph{by} Gunther Cornelissen (Utrecht)}\\[1mm]
{\large \emph{at the} Arbeitstagung 2011  } 

\hfill To Don Zagier, on his 60th birthday\\

\textbf{Counting points.} Let $X$ denote a smooth projective curve over a finite field $k=\F_q$. Is $X$ determined (up to isomorphism) from counting its points over finite extensions of $k$, i.e., by the numbers $N_n:= |X(\F_{q^n})|, $
i.e., by knowing its \emph{zeta function} 
$$ \zeta_X(s):= \exp \left( \sum_{n \geq 1} N_n \frac{q^{-sn}}{n} \right) \mbox{ ?}$$
The answer is \textbf{no} in general. Tate (1966) and Turner (1978) proved that for two curves $X, Y$ over $k$, the equality $\zeta_X=\zeta_Y$ is equivalent to their respective Jacobians $\mathrm{Jac}(X) \sim \mathrm{Jac}(Y)$ being $k$-isogenous. The following example of E.\ Howe from 1996 illustrates this phenomenon: let $ X_{\pm} \colon y^2 = x^5\pm x^3+x^2-x-1$ over $\F_3$. Then 
$$ \zeta_{X_{\pm}} = \frac{1-T+T^2-3T^3+9T^4}{(1-T)(1-3T)} \mbox{ with } T=q^{-s}, $$
and here are the first few point counts (for this occasion done independently in Sage):
\begin{center}
\begin{tabular}{llllllll}
$n$ & 1 & 2 & 3 & 4 & 5 & 6 & 7 \dots \\ \hline
$N_n$ & 3 & 11 & 21 & 107 & 288 & 719 & 2271 \dots 
\end{tabular}
\end{center}
Can we remedy this? 

\textbf{Number fields.} Now consider the same problem for a number field $K$, with its Dedekind zeta function 
$$ \zeta_K(s) := \sum_{0 \neq \mathfrak{a}} \frac{1}{N(\mathfrak{a})^s}, $$
where the sum runs over all non-zero ideals $\mathfrak{a}$ of the ring of integers of $K$. 
Knowing $\zeta_K$ is the same as knowing $f(\mathfrak{p}|p)$ for all prime ideals $\mathfrak{p}$. A Theorem of Mih\'aly Bauer (1903) says that {if $K,L$ are two number fields} \emph{that are Galois over $\Q$}, {then $K \cong L$ is equivalent to $\zeta_K = \zeta_L$}. However, a result of Ga{\ss}mann from 1926 says that in general, there do exist non-isomorphic number fields $K,L$ with $\zeta_K=\zeta_L$. Actually, he proves that $\zeta_K=\zeta_L$ is equivalent to the following statement: fix a common extension $N$ of $K$ and $L$ that is Galois over $\Q$ with Galois group $G$, and let $H_K$ and $H_L$ denote the Galois groups of $N/K$ and $N/L$, respectively. Then $\zeta_K=\zeta_L$ if and only if each $G$-conjugacy class intersects $H_K$ and $H_L$ in the same number of elements. A result from Perlis from 1977 says that the first example with $\zeta_K=\zeta_L$ but $K \not \cong L$ occurs in degree $7$ over $\Q$, and an example is given by $K=\Q(\alpha), L=\Q(\beta)$ with 
$$ \alpha^7-7\alpha+3=0 \mbox{ and } \beta^7+14 \beta^4 - 42 \beta^2 - 21 \beta + 9 = 0. $$ Can we remedy this? 

\textbf{Historical aside: internal/external = failure/success.} Here are some further attempts at finding objects that determine isomorphism of number fields $K$ and $L$: an \emph{isomorphism of adele rings} $ \mathbf{A}_K  \cong \mathbf{A}_L $ is stronger than equality of zeta functions (strictly stronger for number fields, equivalent for function fields), but still does not imply field isomorphism (Komatsu, 1976); an example is $K=\Q(\sqrt[8]{18}) $and $L=\Q(\sqrt[8]{288})$. An isomorphism of abelian Galois groups $G_K^{\ab} \cong G_L^{\ab}$ is not enough, either: Kubota determined the isomorphism type of $G_K^{\ab}$ (its \emph{Ulm invariants}) in terms of $K$, and Onabe (1976) gave explicit examples, such as $G_{\Q(\sqrt{-2})}^\ab \cong G_{\Q(\sqrt{-3})}^\ab$. At the other side of the spectrum, an isomorphism of absolute Galois groups $G_K \cong G_L$ does imply that $K \cong L$! This is due to Neukirch (1969) when $K,L$ are Galois over $\Q$ and Uchida (1976) in general. This last theorem is the first manifestation of what Grothendieck called \textbf{an}abelian theorems. We conclude that the objects listed above, that are \emph{internal} to a number field $K$ (i.e., can be described in terms of ideals of $K$), such as $\zeta_K, \mathbf{A}_K$ or $G_K^{\ab}$ (which is internal by class field theory), lead to \emph{failure}, whereas a mysterious object $G_K$, that is \emph{external} to $K$ (described in terms of extensions of $K$, or via the Langlands program in terms of automorphic forms), leads to \emph{success} \dots Can we do better, and have internal success? 

\textbf{Method: class field theory as (noncommutative) dynamical system.} 
Let $J_K$ denote the group of fractional ideals of $K$, $J_K^+$ the semigroup of integral ideals of $K$, $\vartheta_K \colon \mathbf{A}_K^* \rightarrow G_K^{\ab}$ the Artin reciprocity map and $\hat {\mathcal{O}}_K$ the integral finite adeles of $K$. Choose a section $s$ of the natural map $\mathbf{A}_{K,f}^*\rightarrow J_{K} \colon (x_{\p})_{\p} \mapsto \displaystyle\prod \p^{v_{\p}(x_{\p})}.$

These objects were used by Ha and Paugam in 2005 to construct a dynamical system associated to $K$  (for $K=\Q$, this is the famous Bost-Connes system), as follows: we make a \emph{topological space} 
$$ X_K = G_K^{\ab} \times_{\hat{\mathcal{O}}_K^*} \hat {\mathcal{O}}_K, $$
consisting of classes $[(\gamma,\rho)]$ for $\gamma \in G_{K}^{\ab}$ and $\rho \in \hat{\mathcal{O}}_{K}$, defined by the equivalence $$(\gamma,\rho) \sim (\vartheta_{K}({u}^{-1}) \cdot \gamma, u \rho) \mbox{ for all } u \in \hat{\mathcal{O}}_{K}^*.$$
Then we consider the \emph{action} of $\mathfrak{n} \in J_K^+$ on $X_K$ given by 
$$ \mathfrak{n} \ast [(\gamma,\rho)] := [(\vartheta_{K}(s(\mathfrak{n}))^{-1} \gamma, s(\mathfrak{n}) \rho)].$$ In this way, we get a dynamical system $(X_K, J_K^+)$. 

\textbf{Main Theorem.} (C-Matilde Marcolli,  arxiv:1009.0736) \emph{ For two number fields $K$ and $L$, an isomorphism $K \cong L$ is equivalent to a norm-preserving isomorphism of dynamical systems $(X_K,J_K^+) \cong (X_L,J_L^+).$} 

By \emph{isomorphism of dynamical systems}, we mean a homeomorphism $\Phi \colon X_K \isomto X_L$ and a group homomorphism $\varphi \colon J_K^+ \isomto J_L^+$ such that $\Phi(\mathfrak{n} \ast x ) = \varphi(\mathfrak{n}) \ast \Phi(x)$ for all $x \in X_K$ and $\mathfrak{n} \in J_K^+$; and \emph{norm-preserving} means that $N_L(\varphi(\mathfrak{n}))=N_K(\mathfrak{n})$ for all $\mathfrak{n} \in J_K^+$. 

In a sense, this theorem shows that a \emph{suitable combination of failure} ($\zeta_K$, which will be the partition function of the system, $G_K^{\ab}$ and $\mathbf{A}_K$, which occur in the system) \emph{may lead to success.} It gives an ``internal'' description of the isomorphism type of a number field. It also holds in a function field, with a slightly different, easier proof. 

The proof is really to ``hit the dynamical system with a hammer until enough isomorphic objects jump out''. 

\textbf{Reformulation using Quantum Statistical Mechanics.} There is a way to reformulate the main theorem by encoding the dynamics in Banach algebra language. We set $A_K:= C(X_K) \rtimes J_K^+$ to be the semigroup crossed product $C^*$-algebra corresponding to the dynamical system. Physically, it corresponds to the \emph{algebra of observables}. If we let $\mu_{\mathfrak{n}}$ and $\mu^*_{\mathfrak{n}}$ denote the partial isometries of the algebra corresponding to $\mathfrak{n} \in J_K^+$, then we also need the non-involutive subalgebra $A_K^{\dagger}$ of $A_K$ generated by $C(X)$ and $\langle \mu_{\mathfrak{n}} \rangle_{\mathfrak{n} \in J_K^+}$ (but not the $\mu_{\mathfrak{n}}^*$). We also consider a one-parameter subgroup of automorphisms of $A_K$, denoted $\sigma_K \colon \mathbf{R} \hookrightarrow \mathrm{Aut}(A_K)$, defined by $\sigma_K(t)(f)=f$ and  $\sigma_K(t)(\mu_{\mathfrak{n}}) = N_K(\mathfrak{n})^{it} \mu_{\mathfrak{n}}$. The algebra with this so-called \emph{time evolution} is an abstract \emph{quantum statistical mechanical system.} A slightly stronger statement than the main theorem is the following: $K \cong L$ is equivalent to an isomorphism of $(A_K,\sigma_K) \isomto (A_L,\sigma_L)$ that maps $A_K^{\dagger}$ to $A_L^{\dagger}$.

From the main theorem, we can deduce our answer to the problems outlined before: 

\textbf{Theorem.} \emph{If $K$ and $L$ are global fields (number fields, or function fields of curves over finite fields), then $K \cong L$ (which, in the case of function fields of curves is equivalent to isomorphism of the curves) is equivalent to the existence of an isomorphism}
$ \psi \colon G_{K}^{\ab} \isomto G_{L}^{\ab}, $
\emph{such that \textbf{all} abelian $L$-series match:}
$ L_K(\chi) = L_L((\psi^{-1})^*\chi) \mbox{\  \ \emph{for all} \ } \chi \in \mathrm{Hom}(G_{K}^{\ab},S^1). $

We discovered this theorem because $L$-series occur as evaluations of low temperature equilibrium states of the system at particular test functions related to the character. Our proof of this theorem is to deduce from $L$-series equality an isomorphism of dynamical systems, which basically boils down to a bit of character theory, and then using the main theorem. In the meanwhile, Bart de Smit has discovered a purely number theoretical proof of the theorem for $L$-series for number fields, and has actually proven something much stronger: for every number field $K$, there is a character of order $3$, such that $L_K(\chi) \neq L_{K'}(\chi')$ for \emph{every} number field $K' \not \cong K$ and character for $G_{K'}^{\ab}$. This proof does not seem to transfer readily to function fields. 

Final remark: the theorem is not really an analytic statement. It suffices to have equality of $L$-series at sufficiently large integers. Hence the theorem also holds with $p$-adic $L$-functions. One may read it as an equivalence of rank-one motives over $K$ and $L$.

\textbf{An analog in Riemannian geometry.} The isospectrality problem has a long history, that can be traced back at least to the Wolfskehl lecture of the dutch physicist Lorentz in G\"ottingen in 1910, where he asked whether the spectrum of the Laplacian on a domain (with suitable boundary conditions) determines the volume. He refers to the Leiden PhD thesis of Johanna Reudler, that very cleverly computes several convincing examples (published in 1912). Hermann Weyl proved the general case in 1911, and much later Mark Kac popularized the question whether the entire shape of the region (so up to euclidean transformations) is determined by the spectrum, as ``Can you hear the shape of a drum?''(this formulation is due to Bers, the problem was originally posed by Bochner). The first counterexample was the construction of two non-isometric Riemannian manifolds with the same spectrum by Milnor, based on Witt's theory of quadratic forms. Then even came non-homeomorphic isospectral manifolds in the work of Ikeda (lens spaces) and Vign\'eras (3-manifolds). 

Let $(X,g)$ denote a closed Riemannian manifold with Laplace operator $\Delta_X$.  The question whether or not the spectrum (with multiplicities) determines the isometry type of $X$ is the same as that whether or not the \emph{spectral zeta function} $$ \zeta_X(s)=\sum_{\lambda \neq 0} \frac{1}{\lambda^s}  = \tr(\Delta_X^{-s})$$
(sum over the non-zero eigenvalues of the Laplace operator, with multiplicities) does so. Can we do better? This time, our ``remedy'' is the following: 
for $a \in C(X)$, set $\zeta_{X,a}(s) = \tr(a\Delta_X^{-s})$, and for $a \in W(X)$ (Lipschitz functions) set $\tilde\zeta_{X,a}=\tr(a[\Delta_X,a]\Delta_X^{-s})$. Then: 

\textbf{Theorem.} (C-Jan Willem de Jong; arXiv:1007.0907) \emph{Let $X$ and $Y$ denote two closed RIemannian manifolds, and $\varphi \colon X \rightarrow Y$ a $C^1$-bijective map. Then $\varphi$ being an isometry is equivalent to the following two properties holding simultaneously    \begin{enumerate}
\item[\textup{(a)}] $ \zeta_{Y,a_0} = \zeta_{X,\varphi^*(a_0)}$ for all $a_0 \in C(Y)$, and 
\item[\textup{(b)}] $ \tilde \zeta_{Y,a_1} = \tilde \zeta_{X,\varphi^*(a_1)}$ for all $a_1 \in W(Y)$.
\end{enumerate}
  }

The proof is a rather formal computation with residues. Various analytically more challenging amplifications are possible, for example, condition (a) alone suffices when the spectrum is simple (which is the generic case by a result of Uhlenbeck). In the above theorem, one can also restrict to a countable dense subset of functions, and to sufficiently large integral values of the zeta functions, so the characterisation is really by countably many values.

\textbf{Lengths of maps.} One may now define the \emph{length of a map} $\varphi \colon X \rightarrow Y$ as the ``distance between the (meromorphic) zeta functions that occur in the theorem''. The usual distance of meromorphic functions doesn't quite work, but the following does: The length $\ell(\varphi)$ of $\varphi$ of Riemannian manifolds of dimension $n$ is 
$$ \ell(\varphi):= \mathop{\sup_{a_0 \in C(Y,\mathbf{R}_{\geq 0})-\{0\}}}_{a_1\in W^1(Y)-\mathbf{R}}\ \sup_{n \leq s \leq n+1} \! \! \! \max \, \{ | \log \left| \frac{\zeta_{X,a^*_0}(s)}{\zeta_{Y,a_0}(s)} \right| |, | \log \left| \frac{\tilde\zeta_{X,a^*_1}(s)}{\tilde\zeta_{Y,a_1}(s)} \right| | \}. $$ 
This then satisfies $\ell(\varphi)=0$ if and only if $\varphi$ is an isometry, and $\ell(\psi \circ \varphi) \leq \ell(\psi) + \ell(\varphi).$
One can also show that $$d(X,Y):=\max \{ \inf_{C^1(X \stackrel{\varphi}{\rightarrow} Y)} \ell(\varphi),+\infty\}$$
defines an extended metric between isometry classes of Riemannian manifolds. 

As an example especially for Don Zagier, we bound the distance $d$ between two tori, corresponding to $i$ and $\rho = (1+\sqrt{-3})/2$ in the upper half plane. This will satisfy 
$$ e^d \leq \frac{E(i,2)}{E(\rho,2)} = \frac{\zeta_{m^2+n^2}(2)}{\zeta_{m^2-mn+n^2}(2)} = \frac{3 \sqrt{3}}{4} \cdot \frac{D(i)}{D(\rho)} = 1.17235730884473\dots, $$
where $E$ is an Eisenstein series, $\zeta_Q$ (with $Q$ a binary quadratic form) is the Epstein zeta function, and $D$ is the Bloch-Wigner dilogarithm function.

\textbf{Pluralizing zeta.} ZETA counts things (points, ideals, geodesics, spectra, \dots) --- it is beautiful, but sometimes lonely, it can fail as an \emph{individual.} But it will be happy and succeed as \emph{part of a family} of ZETAS. 

\end{document}